\documentclass[11pt]{article}
\usepackage{amsmath,amssymb}

\title{Hardness and Algorithms for Rainbow Connection}

\author{
Sourav Chakraborty
\thanks{Department of Computer Science,
University of Chicago, Chicago, IL-60637 USA. Email: \hbox{sourav@cs.uchicago.edu}}
\and
Eldar Fischer
\thanks{Department of Computer Science, Technion, Haifa 32000, Israel. Email:
\hbox{eldar@cs.technion.ac.il}}
\and
Arie Matsliah
\thanks{Department of Computer Science, Technion, Haifa 32000, Israel. Email:
\hbox{ariem@cs.technion.ac.il}}
\and
Raphael Yuster
\thanks{Department of Mathematics, University of Haifa, Haifa 31905, Israel. Email:
\hbox{raphy@math.haifa.ac.il}}
}

\date{}

\newtheorem{theorem}{Theorem}[section]
\newtheorem{lemma}[theorem]{Lemma}
\newtheorem{corollary}[theorem]{Corollary}

\newtheorem{obs}[theorem]{Observation}
\newtheorem{definition}[theorem]{Definition}

\newcommand{\qed}{\rule{2mm}{2mm}}
\newenvironment{proof}{\par\noindent{\bf Proof.}\,}{  $\qed$}

\newcommand{\low}{{\bf low}}
\newcommand{\high}{{\bf high}}
\newcommand{\ov}{\overline{v}}
\newcommand{\diam}{\mathrm{diam}}

\newcommand{\neweps}{\epsilon}

\begin{document}
\maketitle
\setcounter{page}{1}
\begin{abstract}
An edge-colored graph $G$ is {\em rainbow connected}
if any two vertices are connected by a path whose edges have distinct colors.
The {\em rainbow connection} of a connected graph $G$, denoted $rc(G)$,
is the smallest number of colors that are needed in order to make $G$ rainbow
connected. In the first result of this paper we prove that computing $rc(G)$ is NP-Hard
solving an open problem from \cite{Ca-Yu}.
In fact, we prove that it is already NP-Complete to decide if $rc(G)=2$, and
also  that it is NP-Complete to decide
whether a given edge-colored (with an unbounded number of colors) graph is rainbow connected.
On the positive side, we prove that for every $\epsilon >0$, a connected graph with
minimum degree
at least $\epsilon n$ has {\em bounded} rainbow connection, where the bound depends only on
$\epsilon$, and a corresponding coloring can be constructed in polynomial time.
Additional non-trivial upper bounds, as well as open problems and
conjectures are also presented.
\end{abstract}

\section{Introduction}

Connectivity is perhaps the most fundamental graph-theoretic property, both in the
combinatorial sense and the algorithmic sense.
There are many ways to strengthen the connectivity property,
such as requiring hamiltonicity, $k$-connectivity, imposing bounds on the diameter,
requiring the existence of edge-disjoint spanning trees, and so on.

An interesting way to quantitavely strengthen the connectivity
requirement was recently introduced by Chartrand et al. in
\cite{Ch-Zh}. An edge-colored graph $G$ is {\em rainbow connected}
if any two vertices are connected by a path whose edges have
distinct colors. Clearly, if a graph is rainbow connected, then it
is also connected. Conversely, any connected graph has a trivial
edge coloring that makes it rainbow connected; just color each edge
with a distinct color. Thus, one can properly define the {\em
rainbow connection} of a connected graph $G$, denoted $rc(G)$, as
the smallest number of colors that are needed in order to make $G$
rainbow connected. An easy observation is that if $G$ is connected and has $n$
vertices then $rc(G) \le n-1$, since one may color the edges of a
given spanning tree with distinct colors. We note also the trivial
fact that $rc(G)=1$ if and only if $G$ is a clique, the (almost)
trivial fact that $rc(G)=n-1$ if and only if $G$ is a tree, and the
easy observation that a cycle with $k > 3$ vertices has rainbow
connection $\lceil k/2 \rceil$. Also notice that, clearly, $rc(G)
\ge \diam(G)$ where $\diam(G)$ denotes the diameter of $G$.

Chartrand et al.\ computed the rainbow connection of several graph
classes including complete multipartite graphs \cite{Ch-Zh}. Caro et
al.\ \cite{Ca-Yu} considered the extremal graph-theoretic aspects of
rainbow connection. They proved that if $G$ is a connected graph
with $n$ vertices and with minimum degree $3$ then $rc(G) < 5n/6$,
and if the minimum degree is $\delta$ then $rc(G) \le \frac{\ln
\delta}{\delta}n(1+f(\delta))$ where $f(\delta)$ tends to zero as
$\delta$ increases. They also determine the threshold function for a
random graph $G(n,p(n))$ to have $rc(G)=2$. In their paper, they
conjecture that computing $rc(G)$ is an NP-Hard problem, as well as
conjecture that even deciding whether a graph has $rc(G)=2$ in
NP-Complete.

In this paper we address the computational aspects of rainbow connection.
Our first set of results solve, and extend, the complexity conjectures from
\cite{Ca-Yu}. Indeed, it turns out that deciding whether $rc(G)=2$ is
an NP-Complete problem. Our proof is by a series of reductions, where on the way
it is shown that $2$-rainbow-colorability is computationally equivalent to the
seemingly harder question of deciding the existence of a $2$-edge-coloring that
is required to rainbow-connect only vertex pairs from a prescribed set.
\begin{theorem}
\label{t-rc2npc}
Given a graph $G$, deciding if $rc(G)=2$ is NP-Complete.
In particular, computing $rc(G)$ is NP-Hard.
\end{theorem}
Suppose we are given an edge coloring of the graph. Is it then easier to verify whether
the colored graph is rainbow connected? Clearly, if the number of colors in constant then
this problem becomes easy. However, if the coloring is arbitrary, the problem becomes
NP-Complete:
\begin{theorem}
\label{t-rcnpc}
The following problem is NP-Complete: Given an edge-colored graph $G$, check whether the
given coloring makes $G$ rainbow connected.
\end{theorem}
For the proof of Theorem \ref{t-rcnpc}, we first show that the
$s-t$ version of the problem is NP-Complete. That is, given two vertices $s$ and $t$ of an
edge-colored graph, decide whether there is a rainbow path connecting them.

We now turn to positive algorithmic results. Our main positive result is that connected
$n$-vertex graphs with minimum degree $\Theta(n)$
have {\em bounded} rainbow connection. More formally, we prove:
\begin{theorem}\label{thm:main}
For every $\epsilon > 0$ there is a constant $C=C(\epsilon)$ such that
if $G$ is a connected graph with $n$ vertices and minimum degree
at least $\epsilon n$, then $rc(G) \le C$.
Furthermore, there is a polynomial time algorithm that constructs
a corresponding coloring for a fixed $\epsilon$.
\end{theorem}
The proof of Theorem \ref{thm:main} is based upon a modified
degree-form version of Szemer\'edi's Regularity Lemma that we prove
and that may be useful in other applications. From our algorithm it is also
not hard to find a probabilistic polynomial time algorithm for finding this
coloring with high probability (using on the way the algorithmic version
of the Regularity Lemma from \cite{Al-Yu} or \cite{Fi-Sh}).

We note that connected graphs with minimum degree $\epsilon n$ have bounded diameter, but
the latter property by itself does {\em not} guarantee bounded rainbow connection. As an
extreme example, a star with $n$ vertices has diameter $2$ but its rainbow connection is
$n-1$.
The following theorem asserts however that having diameter $2$ and only logarithmic minimum degree
suffices to guarantee rainbow connection $3$.
\begin{theorem}
\label{t-8logn}
If $G$ is an $n$-vertex graph with diameter $2$ and minimum degree at least $8\log n$ then
$rc(G) \le 3$. Furthermore, such a coloring is given with high probability by a uniformly
random $3$-edge-coloring of the graph $G$, and can also be found by a polynomial time
deterministic algorithm.
\end{theorem}
Since a graph with minimum degree $n/2$ is connected and has diameter $2$, we have as an
immediate corollary:
\begin{corollary}
If $G$ is an $n$-vertex graph with minimum degree at least $n/2$ then
$rc(G) \le 3$.
\end{corollary}

The rest of this paper is organized as follows. The next section
contains the hardness results, including the proofs of Theorem
\ref{t-rc2npc} and Theorem \ref{t-rcnpc}. Section
\ref{sec:alg} contains the proof of Theorem \ref{thm:main} and the proof of
Theorem \ref{t-8logn}. At the end of the proof of
each of the above theorems we explain how the algorithm can be
derived -- this mostly consists of using the conditional expectation
method to derandomize the probabilistic parts of the proofs. The final Section \ref{sec:concl} contains some open
problems and conjectures.

\section{Hardness results} \label{sec:hard}

We first give an outline of our proof of Theorem \ref{t-rc2npc}.
We begin by showing the computational equivalence of the problem of rainbow connection $2$, that asks
for a red-blue edge coloring in which {\em all} vertex pairs have a rainbow path connecting
them, to the problem of {\em subset rainbow connection $2$}, asking for a
red-blue coloring in which every pair of vertices in a {\em given subset} of pairs has a
rainbow path connecting them. This is proved in Lemma \ref{l-21} below.

In the second step, we reduce the problem of {\em extending to rainbow connection $2$},
asking whether a given partial red-blue coloring can be completed to a obtain a rainbow
connected graph,
to the subset rainbow connection $2$ problem.
This is proved in Lemma \ref{l-22} below.

Finally, the proof of Theorem \ref{t-rc2npc} is completed by reducing $3$-SAT to the
problem of {\em extending to rainbow connection $2$}.

\begin{lemma}
\label{l-21}
The following problems are polynomially equivalent:
\begin{enumerate}
\item
Given a graph $G$ decide whether $rc(G)=2$.
\item
Given a graph $G$ and a set of pairs $P \subseteq V(G) \times V(G)$, decide whether there
is an edge coloring of $G$ with $2$ colors such that all pairs $(u,v) \in P$ are rainbow
connected.
\end{enumerate}
\end{lemma}
\begin{proof}
It is enough to describe a reduction from Problem 2 to Problem 1. Given a graph $G = (V,E)$
and a set of pairs $P \subseteq V \times V$, we construct a graph
$G' = (V',E')$ as follows.

For every vertex $v \in V$ we introduce a new vertex $x_v$, and for
every pair $(u,v) \in (V \times V) \setminus P$ we introduce a new
vertex $x_{(u,v)}$. We set $$V' = V \cup \{x_v:v \in V\} \cup
\{x_{(u,v)}:(u,v) \in (V \times V)\setminus P\}$$ and $$E' = E \cup
\Big\{ \{v,x_v\}:v \in V \Big\} \cup \Big\{
\{u,x_{(u,v)}\},\{v,x_{(u,v)}\}:(u,v) \in (V \times V)\setminus
P\Big\} \cup \Big\{ \{x,x'\}: x,x' \in V' \setminus V\Big\}.$$

It remains to verify that $G'$ is $2$-rainbow connected if and
only if there is an edge coloring of $G$ with $2$ colors such that
all pairs $(u,v) \in P$ are rainbow connected. In one direction, notice that when $G$
is considered as a subgraph of $G'$, no pair of vertices of $G$ that appear in $P$ has
a path of length two in $G'$ that is not fully contained in $G$. Hence, if $G$ is not colorable
in a way connecting the pairs in $P$, the graph $G'$ is not $2$-rainbow-connected.

In the other direction, assume that
$\chi:E\to\{\mathrm{red},\mathrm{blue}\}$ is a coloring of $G$ that
rainbow-connects the pairs in $P$. To extend it to a
rainbow-coloring $\chi':E'\to\{\mathrm{red},\mathrm{blue}\}$, define
$\chi'(\{v,x_v\})=\mathrm{blue}$ for all $v\in V$,
$\chi'(\{u,x_{(u,v)}\})=\mathrm{blue}$ and
$\chi'(\{v,x_{(u,v)}\})=\mathrm{red}$ for all $(u,v)\not\in P$ (note
that we treat $P$ as a set of ordered pairs -- an unordered pair can
be represented by putting both orderings in $P$), and finally
$\chi'(\{x,x'\})=\mathrm{red}$ for all $x,x'\in V'\setminus V$. One
can see that $\chi'$ is indeed a valid $2$-rainbow-coloring of $G'$,
concluding the proof.
\end{proof}

\begin{lemma}
\label{l-22}
The first problem defined below is polynomially reducible to the second one:
\begin{enumerate}
\item
Given a graph $G=(V,E)$ and a partial $2$-edge-coloring $\hat{\chi}: \hat{E} \rightarrow
\{0,1\}$ for $\hat{E}\subset E$, decide whether
$\hat{\chi}$ can be extended to a complete $2$ edge-coloring $\chi:E
\rightarrow \{0,1\}$ that makes $G$ rainbow connected.
\item
Given a graph $G$ and a set of pairs $P \subseteq V(G) \times V(G)$ decide whether there
is an edge coloring of $G$ with $2$ colors such that all pairs $(u,v) \in P$ are rainbow
connected.
\end{enumerate}
\end{lemma}
\begin{proof}
Since the identity of the colors does not matter, it is more
convenient that instead of a coloring $\chi :E \rightarrow \{0,1\}$
we consider the corresponding partition $\pi_{\chi} = (E_1,E_2)$ of
$E$. Similarly, in the case of a partial coloring $\hat{\chi}$, the
pair $\pi_{\hat{\chi}} = (\hat{E_1},\hat{E_2})$ will contain the
corresponding disjoint subsets of $E$ (which may not cover $E$).

Now, given such a partial coloring $\hat{\chi}$ we extend the
original graph $G=(V,E)$ to a graph $G'=(V',E')$, and define a set
$P$ of pairs of vertices such that for the resulting graph the
answer for Problem 2 is ``yes'' if and only if the answer for
Problem 1 for the original graph is ``yes''.

Let $\ell:V \rightarrow [|V|]$ be arbitrary linear ordering of the
vertices, and let $\high:E \rightarrow V$ be a mapping that maps an
edge $e= \{u,v\}$ to $u$ if $\ell(u) > \ell(v)$, and to $v$
otherwise. Similarly, let $\low:E \rightarrow V$ be a mapping that
maps an edge $e= \{u,v\}$ to $u$ if $\ell(u) < \ell(v)$, and to $v$
otherwise.

We construct $G'$ as follows. We add $3 +|\hat{E_1}| + |\hat{E_2}|$
new vertices
$$\{b_1,c,b_2\} \cup \{c_{e} :e \in (\hat{E_1} \cup \hat{E_2})\}$$ and add the
edges $$\Big\{ \{b_1,c\},\{c,b_2\} \Big\} \cup \Big\{\{b_i,c_e\}: i
\in \{1,2\},\ e \in \hat{E}_i \Big\} \cup \Big\{\{c_e,\low(e)\}: e
\in (\hat{E}_1 \cup \hat{E}_2) \Big\}.$$ Now we define the set $P$
of pairs of vertices that have to be $2$-rainbow connected:
$$P = \Big\{ \{b_1,b_2\} \Big\} \cup \Big\{ \{u,v\} : u,v \in V \Big\} \cup \Big\{ \{c,c_e\} : e \in (\hat{E_1} \cup
\hat{E_2})\Big\}
$$ $$ \cup \Big\{ \{b_{i}, \low(e)\} : i \in \{1,2\},\ e \in \hat{E_i}\Big\}
\cup \Big\{ \{c_e, \high(e)\} : e \in (\hat{E_1} \cup
\hat{E_2})\Big\}.
$$

Now, if there is a $2$-rainbow-coloring $\pi_\chi=(E_1,E_2)$ of $G$
which extends $\pi_{\hat{\chi}}=(\hat{E_1},\hat{E_2})$, then we
color $G'$ as follows. $E$ is colored as mandated by $\chi$, that
is, an edge is colored red if it is in $E_1$ and otherwise it is
colored blue. Now $\{b_1,c\}$, $\{b_2,c_e\}$ for $e\in\hat{E_2}$ and
$\{c_e,\mathbf{low}(e)\}$ for $e\in\hat{E_1}$ are all colored blue,
and $\{b_2,c\}$, $\{b_1,c_e\}$ for $e\in\hat{E_1}$ and
$\{c_e,\mathbf{low}(e)\}$ for $e\in\hat{E_2}$ are all colored red.
One can see that this coloring indeed rainbow-connects all the pairs
in $P$.

On the other hand, any $2$-edge-coloring of $G'$ that connects the
pairs in $P$ clearly contains a $2$-rainbow-coloring of $G$, because
$P$ contains all vertex pairs of $G$ when considered as an induced
subgraph of $G'$, and also $G'$ contains no path of length $2$
between vertices of $G$ that is not contained in $G$. Also, such a
coloring would have to color $\{c,b_1\}$ and $\{c,b_2\}$
differently. It would also have to color every $\{b_i,c_e\}$ in a
color different from that of $\{c,b_i\}$, and would hence color
$\{c_e,\mathbf{low}(e)\}$ in a color identical to that of
$\{c,b_i\}$, because it has to be the color different from that of
$\{b_i,c_e\}$. Finally, every $e\in\hat{E_i}$ would have to be
colored with the color different from that of $\{c,b_i\}$ so as to
rainbow-connect $\mathbf{high}(e)$ and $c_e$. This means that the
coloring of $G'$ not only provides a $2$-rainbow-coloring of $G$,
but that it also conforms to the original partial coloring
$\hat{\chi}$.
\end{proof}

\medskip \noindent
{\bf Proof of Theorem \ref{t-rc2npc}.}\,
We show that Problem 1 of Lemma \ref{l-22}
is NP-hard, and then deduce that $2$-rainbow-colorability is NP-Complete
by applying Lemma \ref{l-21} and Lemma \ref{l-22} while observing that it clearly belongs to
NP.

We reduce $3$-SAT to Problem 1 of Lemma \ref{l-22}.
Given a 3CNF formula $\phi = \bigwedge_{i=1}^m c_i$ over variables
$x_1,x_2,\ldots ,x_n$, we construct a graph $G_\phi$ and a partial $2$-edge coloring
$\chi':E(G_\phi) \rightarrow \{0,1\}$ such that there is an
extension $\chi$ of $\chi'$ that makes $G_\phi$ rainbow connected if and only if $\phi$ is
satisfiable.

We define $G_\phi$ as follows:
$$ V(G_\phi) = \{c_i : i \in [m]\} \cup \{x_i : i \in [n] \} \cup
\{a\}$$
$$ E(G_\phi) = \Big\{\{c_i,x_j\} : x_j \in c_i \mathrm{\ in\ }\phi \Big\} \cup
\Big\{\{x_i,a\} : i \in [n] \Big\} \cup
\Big\{\{c_i,c_j\} : i,j \in [m] \Big\} \cup \Big\{\{x_i,x_j\} : i,j
\in [n] \Big\}$$ and we define the partial coloring $\chi'$ as
follows:
$$\forall_{i,j \in [m]} \chi'(\{c_i,c_j\}) = 0$$
$$\forall_{i,j \in [n]} \chi'(\{x_i,x_j\}) = 0$$
$$\forall_{\{x_i,c_j\} \in E(G_\phi)} \chi'(\{x_i,c_j\}) = 0 \mathrm{\ if\ }x_i\mathrm{\
is\ positive\ in\ }c_j,\ 1 \mathrm{\ otherwise}$$
while all the edges in $\Big\{\{x_i,a\} : i \in [n] \Big\}$ (and only they) are left
uncolored.

Assuming without loss of generality that all variables in $\phi$ appear both as positive
and as negative, one can verify that a $2$-rainbow-coloring
of the uncolored edges corresponds to a satisfying assignment of
$\phi$ and vice versa.
\qed

The proof of Theorem \ref{t-rcnpc} is based upon the proof of the following theorem.
\begin{theorem}
\label{t-st} The following problem is NP-complete: Given an edge
colored graph $G$ and two vertices $s,t$ of $G$, decide whether
there is a rainbow path connecting $s$ and $t$.
\end{theorem}
\begin{proof}
Clearly the problem is in NP. We prove that it is NP-Complete by reducing 3-SAT to it.
Given a 3CNF formula
$\phi = \bigwedge_{i=1}^m c_i$ over variables $x_1,x_2,\ldots,x_n$, we construct a graph
$G_\phi$ with two special vertices $s,t$ and a coloring
$\chi:E(G_\phi) \rightarrow [|E(G_\phi)|]$ such that there is a
rainbow path connecting $s$ and $t$ in $G_\phi$ if and only if $\phi$ is
satisfiable.

We start by constructing an auxiliary graph $G'$ from $\phi$. The
graph $G'$ has $3m + 2$ vertices, that are partitioned into $m + 2$
layers $V_0,V_1,\ldots,V_m,V_{m+1}$, where $V_0 = \{s\}$, $V_{m+1} =
\{t\}$ and for each $i \in [m]$, the layer $V_i$ contains the three
vertices corresponding to the literals of $c_i$ (a clause in
$\phi$). The edges of $G'$ connect between all pairs of vertices
residing in consecutive layers. Formally,
$$
E(G') = \Big\{\{u,v\} : \exists{i \in [m+1]}\mathrm{\ s.t.\ }
u \in V_{i-1} \mathrm{\ and\ } v \in V_i\Big\}.
$$
Intuitively, in our final colored graph $G_\phi$, every rainbow path
from $s$ to $t$ will define a satisfying assignment of $\phi$ in a
way that for every $i \in [m]$, if the rainbow path contains a
vertex $v \in V_i$ then the literal of $c_i$ that corresponds to $v$
is satisfied, and hence $c_i$ is satisfied. Since any path from $s$
to $t$ must contain at least one vertex from every layer $V_i$, this
will yield a satisfying assignment for the whole formula $\phi$. But
we need to make sure that there are no contradictions in this
assignment, that is, no opposite literals are satisfied together.
For this we modify $G'$ by replacing each literal-vertex with a
gadget, and we define an edge coloring for which rainbow paths
yield only consistent assignments.

For every variable $x_j$, $j\in [n]$, let $v_{j_1},v_{j_2},\ldots,v_{j_k}$ be the vertices
of $G'$ corresponding to the positive literal $x_j$, and let
$\ov_{j_1},\ov_{j_2},\ldots,\ov_{j_\ell}$ be the vertices
corresponding to the negative literal $\overline{x}_j$. We can
assume without loss of generality that both $k\geq 1$ and $\ell \geq
1$, since otherwise the formula $\phi$ can be simplified. For every
such variable $x_j$ we also introduce $k \times \ell$ distinct
colors $\alpha^j_{1,1},\ldots,\alpha^j_{k,\ell}$. Next, we transform
the auxiliary graph $G'$ into the final graph $G_\phi$.

For every $a \in [k]$ we replace the vertex $v_{j_a}$ that resides
in layer (say) $V_i$ with $\ell+1$ new vertices
$v_1,v_2,\ldots,v_{\ell+1}$ that form a path in that order. We also
connect all vertices in $V_{i-1}$ to $v_1$ and connect all vertices
in $V_{i+1}$ to $v_{\ell+1}$. For every $b \in [\ell]$, we color the
edge $\{v_b,v_{b+1}\}$ in the new path with the color
$\alpha^j_{a,b}$. Similarly, for every $b \in [\ell]$ we replace the
vertex $\ov_{j_b}$ from layer (say) $V_{i'}$ with $k+1$ new vertices
$\ov_1,\ov_2,\ldots,\ov_{k+1}$ that form a path, and connect all
vertices in $V_{i'-1}$ to $\ov_1$ and all vertices in $V_{i'+1}$ to
$\ov_{k+1}$. For every $a \in [k]$, we color the edge
$\{v_a,v_{a+1}\}$ with $\alpha^j_{a,b}$. All other edges of $G_\phi$
(which were the original edges of $G'$) are colored with fresh
distinct colors.

Clearly, any path from $s$ to $t$ in $G_\phi$ must contain at least
one of the newly built paths in each layer. On the other hand, it
is not hard to verify that any two paths of opposite literals of
the same variable have edges sharing the same color.
\end{proof}
\medskip
\\{\bf Proof of Theorem \ref{t-rcnpc}.}\, We reduce from the problem
in Theorem \ref{t-st}. Given an edge colored graph $G=(V,E)$ with
two special vertices $s$ and $t$, we construct a graph $G'=(V',E')$
and define a coloring $\chi':E' \rightarrow [|E'|]$ of its edges
such that $s$ and $t$ are rainbow connected in $G$ if and only if
the coloring of $G'$ makes $G'$ rainbow connected.

Let $V = \{v_1=s,v_2,\ldots,v_{n}=t\}$ be the vertices of the
original graph $G$. We set $$V' = V \cup \{s',t',b\} \cup
\{s^1,v_2^1,v_2^2,\ldots,v_{n-1}^1,v_{n-1}^2,t^2\}$$ and
$$E' = E \cup
\Big\{\{s',s\},\{t',t\},\{s,s^1\},\{t,t^2\}\Big\} \cup
\Big\{\{b,v_i\}: i \in [n]\Big\} \cup $$ $$\cup
\Big\{\{v_i,v_i^j\}:i \in [n],\ j \in \{1,2\}\Big\} \cup
\Big\{\{v_i^a,v_{j}^{b}:i,j \in [n],\ a,b \in \{1,2\}\Big\}.$$ The
coloring $\chi'$ is defined as follows:

\begin{itemize}
\item all edges $e\in E$ retain the original color, that is $\chi'(e) =
\chi(e)$;
\item the edges $\{t,t'\},\{s,b\}$ and $\Big\{\{v_i,v_i^1\}: i \in [n-1]\Big\}$ are
colored with
a special color $c_1$;
\item the edges $\{s,s'\},\{t,b\}$ and $\Big\{\{v_i,v_i^2\}: i \in [2,n]\Big\}$ are colored
with a special color $c_2$;
\item the edges in $\Big\{\{v_i,b\}: i \in [2,n-1]\Big\}$ are colored
with a special color $c_3$;
\item the edges in $\Big\{\{v_i^a,v_j^b\}: i,j \in [n],\ a,b \in \{1,2\}\Big\}$ are colored
with a special color $c_4$.
\end{itemize}
One can verify that $\chi'$ makes $G'$ rainbow connected if
and only if there was a rainbow path from $s$ to $t$ in $G$. \qed

\section{Upper bounds and algorithms} \label{sec:alg}
The proof of our main Theorem \ref{thm:main} is based upon a
modified degree-form version of Szemer\'edi's Regularity Lemma, that
we prove here and that may be useful in other applications. We begin
by introducing the Regularity Lemma and the already known
degree-form version of it.
\subsection{Regularity Lemma}

The Regularity Lemma of Szemer\'edi \cite{Sz} is one of the most
important results in graph theory and combinatorics, as it
guarantees that every graph has an $\epsilon$-approximation of
constant descriptive size, namely a size that depends only on
$\epsilon$ and not on the size of the graph. This approximation
``breaks'' the graph into a constant number of pseudo-random
bipartite graphs. This is very useful in many applications since
dealing with random-like graphs is much easier than dealing with
arbitrary graphs. In particular, as we shall see, the Regularity
Lemma allows us to prove that graphs with linear minimum degree have
bounded rainbow connection.

We first state the lemma. For two nonempty disjoint vertex sets $A$
and $B$ of a graph $G$, we define $E(A,B)$ to be the set of edges of
$G$ between $A$ and $B$. The {\em edge density} of the pair is
defined by $d(A,B)=|E(A,B)|/(|A||B|)$.

\begin{definition}[$\epsilon$-regular pair]\label{DefSubReg}
A pair $(A,B)$ is {\em $\epsilon$-regular} if for every $A'\subseteq
A$ and $B'\subseteq B$ satisfying $|A'|\geq \epsilon |A|$ and
$|B'|\geq \epsilon |B|$, we have $|d(A',B') - d(A,B)| \leq
\epsilon$.
\end{definition}

An $\epsilon$-regular pair can be thought of as a pseudo-random
bipartite graph in the sense that it behaves almost as we would
expect from a random bipartite graph of the same density.
Intuitively, in a random bipartite graph with edge density $d$, all
large enough sub-pairs should have similar densities.

A partition $V_1,\ldots,V_k$ of the vertex set of a graph is called
an {\em equipartition} if $|V_i|$ and $|V_{j}|$ differ by no more
than $1$ for all $1\leq i < j \leq k$ (so in particular every $V_i$
has one of two possible sizes). The {\em order} of an equipartition
denotes the number of partition classes ($k$ above). An
equipartition $V_1,\ldots,V_k$ of the vertex set of a graph is
called {\em $\epsilon$-regular} if all but at most $\epsilon\binom{k}{2}$ of the pairs $(V_i,V_{j})$ are $\epsilon$-regular.
Szemer\'edi's Regularity Lemma can be formulated as follows.

\begin{lemma}[Regularity Lemma \cite{Sz}]\label{SzReg}
For every $\epsilon > 0$ and positive integer $K$, there exists
$N=N_{\ref{SzReg}}(\epsilon,K)$, such that any graph with $n \geq N$
vertices has an $\epsilon$-regular equipartition of order $k$, where
$K \leq k \leq N$.
\end{lemma}
As mentioned earlier, the following variation of the lemma comes
useful in our context.
\begin{lemma}[Regularity Lemma - degree form \cite{Ko-Si}] \label{deg_form}
For every $\epsilon > 0$ and positive integer $K$ there is
$N=N_{\ref{deg_form}}(\epsilon,K)$ such that given any graph $G =
(V,E)$ with $n > N$ vertices, there is a partition of the vertex-set
$V$ into $k+1$ sets $V_0',V'_1, \ldots, V'_{k}$, and there is a
subgraph $G'$ of $G$ with the following properties:

\begin{enumerate}
\item $K \leq k \leq N$,
\item $s \triangleq |V'_0| \leq \epsilon^5 n$ and all other components $V'_i$, $i \in [k]$ are of size $\ell \triangleq \frac{n-s}{k}$,
\item for all $i\in [k]$, $V'_i$ induces an independent set in $G'$,
\item for all $i,j \in [k]$, the pair $(V'_i,V'_j)$ is $\epsilon^5$-regular in $G'$,
with density either $0$ or at least $\frac{\epsilon}{4}$,
\item for all $v \in V$, $\deg_{G'}(v) > \deg_G(v) - \frac{\epsilon}{3} n$.
\end{enumerate}
\end{lemma}

This form of the lemma (see e.g. \cite{Ko-Si}) can be obtained by
applying the original Regularity Lemma (with a smaller value of
$\epsilon$), and then ``cleaning'' the resulting partition. Namely,
adding to the exceptional set $V'_0$ all components $V_i$ incident
to many irregular pairs, deleting all edges between any other pairs
of clusters that either do not form an $\epsilon$-regular pair or
they do but with density less than $\epsilon$, and finally adding to
$V_0$ also vertices whose degree decreased too much by this deletion
of edges.

\subsection{A modified degree form version of the Regularity Lemma}\label{sec:regularity}
In order to prove that graphs with linear minimum degree have
bounded rainbow connection, we need a special version of
the Regularity Lemma, which is stated next.

\begin{lemma}[Regularity Lemma - new version]\label{lem:deg_new}
For every $\epsilon > 0$ and positive integer $K$ there is
$N=N_{\ref{lem:deg_new}}(\epsilon, K)$ so that the following holds:
If $G = (V,E)$ is a graph with $n > N$ vertices and minimum degree
at least $\epsilon n$ then there is a subgraph $G''$ of $G$, and a
partition of $V$ into $V''_1,\ldots,V''_k$ with the following
properties:
\begin{enumerate}
\item $K \leq k \leq N$,
\item for all $i \in [k]$, $(1-\epsilon)\frac{n}{k} \leq |V''_i| \leq (1 + \epsilon^3)\frac{n}{k}$,
\item for all $i \in [k]$, $V''_i$ induces an independent set in $G''$,
\item for all $i,j \in [k]$, $(V''_i, V''_j)$ is an $\epsilon^3$-regular pair in $G''$,
with density either $0$ or at least $\frac{\epsilon}{16}$,
\item for all $i \in [k]$ and every $v \in V''_i$ there is at least one other class $V''_j$ so that
the number of neighbors of $v$ in $G''$ belonging to $V''_j$ is at
least $\frac{\epsilon}{2}|V''_j|$.
\end{enumerate}
\end{lemma}
\begin{proof} Given $\epsilon > 0$ and $K$, let
us apply the degree-form Regularity Lemma (Lemma \ref{deg_form}),
with the parameters $\neweps$ and $K$, and let $N =
N_{\ref{deg_form}}(\neweps,K)$. Let $V'_0,V'_1,\ldots,V'_k$ be the
partition promised by Lemma \ref{deg_form}, and let $G'$ be the
corresponding subgraph of $G$. Recall the parameters $s \leq
\neweps ^5 n$ and $\ell = \frac{n-s}{k}$ defined in Lemma
\ref{deg_form}.

Consider a mapping $f : V'_0 \rightarrow [k]$ where $f(v) = i$
implies that $v$ has more than $\frac{1}{2} \epsilon \ell$ neighbors
in $G'$ that belong to $V'_i$. Such a mapping clearly exists since
otherwise the degree of $v$ in $G$ would have been at most $(s - 1 +
k \frac{1}{2} \epsilon \ell) + \frac{\neweps}{3} n < \epsilon n$.

Next consider the unique mapping $g : [k] \rightarrow
\mathcal{P}([k])$ that maps an index $i \in [k]$ to the subset of
indices $g(i) \subseteq [k]$ such that $j \in g(i)$ if and only if
$d(V'_i, V'_j) \geq
\neweps / 4$ in $G'$. We notice that the cardinality of $g(i)$ is not
too small. Indeed, if $v$ is any vertex of $V'_i$ then the
cardinality of $g(i)$ is at least $$ \frac{(\deg_{G'}(v)-s)}{\ell} >
\frac{\epsilon }{2}k.$$

Now consider the following process that creates $G''$ from $G'$ by
placing the vertices of $V'_0$ one-by-one in the other vertex
classes. Assume that $V'_0=\{v_1,\ldots,v_s\}$ and denote by $V_i^t$
the extension of $V_i$ after placing $v_t$ somewhere. Initially,
$V_i^0=V'_i$ for $i=1,\ldots,k$, and eventually $V_i^s=V''_i$ for
$i=1,\ldots,k$ forming the desired partition.

At stage $t$, we want to place $v_t \in V'_0$ somewhere. Consider
the subset of indices $g(f(v_t))$. We will choose $j \in g(f(v_t))$
so that $V_j^{t-1}$ has the smallest cardinality (if there are
several candidates, we may choose any of them). We define $V_j^t =
V_j^{t-1} \cup \{v_t\}$ and define $V_i^t=V_i^{t-1}$ for $i \neq j$.
The neighbors of $v_t$ that are kept in $G''$ are all the neighbors
of $v_t$ in $G'$ that belong to $V_x$, where $x \in g(j)$. Notice
that in particular $f(v_t) \in g(j)$.

Having created $G''$ it remains to prove its claimed properties.
Clearly, properties $1$ and $3$ hold. For property $2$, notice that
the final sets $V''_i = V_i^s$ have grown from the initial $V'_i$ by
no more than
$$ \frac{s}{(\frac{\epsilon} {2} k)} < \frac{2}{\epsilon k} \neweps ^5 n = 2\epsilon^4 \frac{n}{k} $$ and
hence $$ \ell = |V'_i| \leq |V''_i| = |V_i^s| \leq \ell +
2\epsilon^4 \frac{n}{k} \leq (1 + \epsilon^3)\ell$$ and property $2$
holds.

For property $4$, Let us first prove the requirement on the density.
Notice that if the original density of $(V'_i,V'_j)$ in $G'$ was
$0$, it remained so also in $G''$. Otherwise, if the density of
$(V'_i,V'_j)$ in $G'$ was at least $\neweps /4$, the fact that
$|V''_i| \leq |V'_i| + 2\epsilon^4 \frac{n}{k} < 2|V'_i|$ implies
that the density of $(V''_i,V''_j)$ is larger than
$\frac{1}{16}\neweps$, no matter where the neighbors of the newly
added vertices reside.

As for the regularity condition in property $4$, notice that if
$(V'_i,V'_j)$ had density $0$ in $G'$, then the same holds in $G''$
and the pair $(V'_i,V'_j)$ is trivially $\epsilon^3$-regular. Now,
if $(V'_i,V'_j)$ had density $\delta \geq \epsilon/4$ in $G'$, there
are precisely $\delta \ell^2$ edges between them. In $G''$ the
number of edges between them can increase by at most $$4 \epsilon^4
\frac{n}{k} (\ell + 2\epsilon^4 \frac{n}{k}) \leq 4 \epsilon^4
\frac{n^2}{k^2}$$ hence the density may have increased from $\delta$
to at most $\delta+4\epsilon^4$. Since $\delta
> \epsilon/4$ and since the pair $(V'_i,V'_j)$ was initially
$\epsilon^5$-regular, it follows that the final pair $(V''_i,V''_j)$
is at least $\epsilon^3$-regular.

For property $5$, this is easily checked to hold for vertices not
originally from $V'_0$ because of property $5$ of $G'$. Property $5$
also holds for vertices of $V'_0$ since when $v_t$ is placed in some
$V'_j$, we know that it has at least $\frac{1}{2} \epsilon \ell$
neighbors in $V'_i$ where $i=f(v_t)$.
\end{proof}
We note that the above a partition as guaranteed by our
modified version of the Regularity Lemma can be found in polynomial
time for a fixed $\epsilon$ (with somewhat worse constants), by
using the exact same methods that were used in \cite{Al-Yu} for
constructing an algorithmic version of the original Regularity
Lemma.

\subsection{Proof of Theorem \ref{thm:main}}

In this section we use our version of the Regularity Lemma to prove
Theorem \ref{thm:main}. First we need some definitions. Given a
graph $G = (V,E)$ and two subsets $V_1,V_2 \subseteq V$, let
$E(V_1,V_2)$ denote the set of edges having one endpoint in $V_1$
and another endpoint in $V_2$. Given a vertex $v$, let $\Gamma(v)$
denote the set of $v$'s neighbors, and for $W \subseteq V$, let
$\Gamma_{W}(v)$ denote the set $W \cap \Gamma(v)$.

For an edge coloring $\chi:E \rightarrow \mathcal{C}$, let
$\pi_\chi$ denote the corresponding partition of $E$ into (at most)
$|\mathcal{C}|$ components. For two edge colorings $\chi$ and
$\chi'$, we say that $\chi'$ is a {\em refinement} of $\chi$ if
$\pi_{\chi'}$ is a refinement of $\pi_\chi$, which is equivalent
to saying that $\chi'(u)=\chi'(v)$ always implies $\chi(u)=\chi(v)$.

\begin{obs} \label{obs:refinement}
Let $\chi$ and $\chi'$ be two edge-colorings of a graph $G$, such
that $\chi'$ is a refinement of $\chi$. For any path $P$ in $G$, if
$P$ is a rainbow path under $\chi$, then $P$ is a rainbow path under
$\chi'$. In particular, if $\chi$ makes $G$ rainbow connected, then so does
$\chi'$. $\qed$
\end{obs}

We define a set of eight distinct colors $\mathcal{C} =
\{a_1,a_2,a_3,a_4,b_1,b_2,b_3,b_4\}$. Given a coloring $\chi:E
\rightarrow \mathcal{C}$ we say that $u,v \in V$ are {\em
$a$-rainbow connected} if there is a rainbow path from $u$ to $v$
using only the colors $a_1,a_2,a_3,a_4$. We similarly define {\em
$b$-rainbow connected} pairs. The following is a central lemma in
the proof of Theorem \ref{thm:main}.

\begin{lemma}\label{lem:ab-rainbow}
For any $\epsilon > 0$, there is
$N=N_{\ref{lem:ab-rainbow}}(\epsilon)$ such that any connected graph
$G = (V,E)$ with $n > N$ vertices and minimum degree at least
$\epsilon n$ satisfies the following. There is a partition $\Pi$ of
$V$ into $k \leq N$ components $V_1,V_2,\ldots,V_k$, and a coloring
$\chi:E \rightarrow \mathcal{C}$ such that for every $i \in [k]$ and
every $u,v \in V_i$, the pair $u,v$ is both $a$-rainbow connected
and $b$-rainbow connected under $\chi$.
\end{lemma}

\begin{proof} Proof given in Section \ref{sec:lempf-ab-rainbow}.
\end{proof}

\

Using Lemma \ref{lem:ab-rainbow} we derive the proof of Theorem
\ref{thm:main}. For a given $\epsilon > 0$, set $N =
N_{\ref{lem:ab-rainbow}}(\epsilon)$ and set $C = \frac{3}{\epsilon}N
+ 8$. Clearly, any connected graph $G = (V,E)$ with $n \leq C$
vertices satisfies $rc(G) \leq C$. So we assume that $n > C \geq N$,
and let $\Pi = V_1,\ldots,V_k$ be the partition of $V$ from Lemma
\ref{lem:ab-rainbow}, while we know that $k \leq N$.

First observe that since the minimal degree of $G$ is $\epsilon n$,
the diameter of $G$ is bounded by $3/\epsilon$. This can be verified
by e.g. by taking an arbitrary vertex $r \in V$ and executing a
$BFS$ algorithm from it. Let $L_1,\ldots,L_t$ be the layers of
vertices in this execution, where $L_i$ are all vertices at distance
$i$ from $r$. Observe that since the minimal degree is at least
$\epsilon n$, the total number of vertices in every three
consecutive layers must be at least $\epsilon n$, thus $t \leq
3/\epsilon$. Since the same claim holds for any $r \in V$, this implies
that $\diam(G) \leq t \leq 3/\epsilon$.

Now let $T = (V_T,E_T)$ be a connected subtree of $G$ on at most $k
\cdot \diam(G) \leq \frac{3}{\epsilon} N$ vertices such that for
every $i \in [k]$, $V_T \cap V_i \ne \emptyset$. Such a subtree must
exist in $G$ since as observed earlier, $\diam(G) \leq 3/\epsilon$.
Let $\chi :E \rightarrow \mathcal{C}$ be the coloring from Lemma
\ref{lem:ab-rainbow}, and let $\mathcal{H} =
\{h_1,h_2,\ldots,h_{|E_T|}\}$ be a set of $|E_T| \leq
\frac{3}{\epsilon} N$ fresh colors. We refine $\chi$ by recoloring
every $e_i \in E(T)$ with color $h_i \in \mathcal{H}$. Let $\chi':E
\rightarrow \Big( \mathcal{C} \cup \mathcal {H} \Big)$ be the
resulting coloring of $G$. The following lemma completes the proof
of Theorem \ref{thm:main}.

\begin{lemma}
The coloring $\chi'$ makes $G$ rainbow connected. Consequently,
$rc(G) \leq |E_T| + 8 \leq C$.
\end{lemma}
\begin{proof}
Let $u,v \in V$ be any pair of $G$'s vertices. If $u$ and $v$ reside
in the same component $V_i$ of the partition $\Pi$, then (by Lemma
\ref{lem:ab-rainbow}) they are connected by a path $P$ of length at
most four, which is a rainbow path under the the original coloring
$\chi$. Since $\chi'$ is a refinement of $\chi$, the path $P$
remains a rainbow path under $\chi'$ as well (see Observation
\ref{obs:refinement}).

Otherwise, let $u \in V_i$ and $v \in V_j$ for $i \neq j$. Let $t_i$
and $t_j$ be vertices of the subtree $T$, residing in $V_i$ and
$V_j$ respectively. By definition of $\chi'$, there is a rainbow
path from $t_i$ to $t_j$ using colors from $\mathcal{H}$. Let $P_t$
denote this path. In addition, by Lemma \ref{lem:ab-rainbow} we know
that for the original coloring $\chi$, there is a rainbow path $P_a$
from $u$ to $t_i$ using colors $a_1,\ldots,a_4$ and there is a
rainbow path $P_b$ from $v$ to $t_j$ using colors $b_1,\ldots,b_4$.
Based on the fact that $\chi'$ is a refinement of $\chi$, it is now
easy to verify that $P_t,P_a$ and $P_b$ can be combined to form a
rainbow path from $u$ to $v$ under $\chi'$.
\end{proof}

\

This concludes the proof of Theorem \ref{thm:main}, apart from the existence
of a polynomial time algorithm for finding this coloring. We note that all
arguments above apart from Lemma \ref{lem:ab-rainbow} admit polynomial
algorithms for finding the corresponding structures. The algorithm for
Lemma \ref{lem:ab-rainbow} will be given with its proof.

\subsection{Proof of Lemma
\ref{lem:ab-rainbow}}\label{sec:lempf-ab-rainbow} First we state
another auxiliary lemma, which is proved in the next section.

\begin{lemma}\label{lem:many}
For every $\epsilon > 0$ there exists $N =
N_{\ref{lem:many}}(\epsilon)$ such that any graph $G = (V,E)$ with
$n > N$ vertices and minimum degree at least $\epsilon n$ satisfies
the following: There exists a partition $\Pi = V_1,\ldots,V_k$ of
$V$ such that for every $i \in [k]$ and every $u,v \in V_i$, the
number of edge disjoint paths of length at most four from $u$ to $v$
is larger than $8^5 \log n$. Moreover, these sets can be found using
a polynomial time algorithm for a fixed $\epsilon$.
\end{lemma}

\begin{proof}{\bf (of Lemma \ref{lem:ab-rainbow})}
First we apply Lemma \ref{lem:many} to get the partition $\Pi$. Now
the proof follows by a simple probabilistic argument. Namely, we
color every edge $e \in E$ by choosing one of the colors in
$\mathcal{C} = \{a_1,\ldots,a_4,b_1,\ldots,b_4\}$ uniformly and independently at
random. Observe that a fixed path $P$ of length at most four is an
$a$-rainbow path with probability at least $8^{-4}$. Similarly, $P$
is a $b$-rainbow path with probability at least $8^{-4}$. So any
fixed pair $u,v \in V_i$ is not both $a$-rainbow-connected and
$b$-rainbow-connected with probability at most $2 (1 -8^{-4})^ {8^5
\log n} < n^{-2}$, and therefore the probability that all such pairs
are both $a$-rainbow connected and $b$-rainbow connected is strictly
positive. Hence the desired coloring must exist.

To find the coloring algorithmically, we note that for every
{\em partial} coloring of the edges of the graph it is easy to
calculate the {\em conditional} probability that the fixed pair
of vertices $u,v$ is not both $a$-rainbow-connected and $b$-rainbow-connected.
Therefor we can calculate the conditional expectation of the number of pairs
that are not so connected for any partial coloring. Now we can
derandomize the random selection of the coloring above by using
the conditional expectation method (cf.\ \cite{AlSp}): In every stage we color one
of the remaining edges in a way that does not increase the conditional
expectation of the number of unconnected pairs. Since this expectation
is smaller than $1$ in the beginning, in the end we will have less
than $1$ unconnected pair, and so all pairs will be connected.
\end{proof}

\subsection{Proof of Lemma \ref{lem:many}}

Given $\epsilon > 0$ let $L = N_{\ref{lem:deg_new}}(\epsilon,1)$
and set $N$ to be the smallest number that satisfies $\epsilon^4
\frac{N}{L} > 8^5 \log N$. Now, given any graph $G = (V,E)$ with
$n > N$ vertices and minimum degree at least $\epsilon n$, we
apply Lemma \ref{lem:deg_new} with parameters $\epsilon$ and $1$.
Let $\Pi = V_1,V_2,\ldots,V_k$ be the partition of $V$ obtained
from Lemma \ref{lem:deg_new}, while as promised, $k \leq L =
N_{\ref{lem:deg_new}}(\epsilon)$.

Fix $i \in [k]$ and $u,v \in V_i$. From Lemma \ref{lem:deg_new} we
know that there is a component $V_a$ such that $u$ has at least
$\frac{\epsilon }{3k}n$ neighbors in $V_a$. Similarly, there is a
component $V_b$ such that $v$ has at least $\frac{\epsilon }{3k}n$
neighbors in $V_b$. Let $\Gamma_{u,a}$ denote the set of $u$'s
neighbors in $V_a$, and similarly, let $\Gamma_{v,b}$ denote $v$'s
neighbors in $V_b$. We assume in this proof that $V_a \neq V_b$,
and at the end it will be clear that the case $V_a = V_b$ can only
benefit.

We say that a set $W_u = \{w_1,\ldots,w_t\} \subseteq V_i$ is {\em
distinctly reachable from $u$} if there are distinct vertices
$w'_1,\ldots,w'_t \in \Gamma_{u,a}$ such that for every $j \in
[t]$, $\{w_j,w'_j\} \in E$. Notice that the collection of pairs
$\{w_j,w'_j\}$ corresponds to a matching in the graph $G$, where
all edges of the matching have one endpoint in $V_i$ and the other
endpoint in $\Gamma_{u,a}$. Similarly, we say that $W_v \subseteq
V_i$ is distinctly reachable from $v$ if there are distinct
vertices $w'_1,\ldots,w'_t \in \Gamma_{v,b}$ such that for every
$j \in [t]$, $\{w_j,w'_j\} \in E$. Observe that it is enough to
prove that there exists a set $W \subseteq V_i$ of
size $\epsilon^4 \frac{N}{L} > 8^5 \log N$ which is
distinctly reachable from both $u$ and $v$.
This will imply the existence of $8^5 \log N$ edge disjoint paths
of length four from $u$ to $v$.

Our first goal is to bound from below the size of the maximal set
$W_u$ as above. Since (by Lemma \ref{lem:deg_new}) $V_a$ and $V_i$
are $\epsilon^3$-regular pairs with density $\geq
\frac{\epsilon}{16}$ and since $\epsilon^3 < \epsilon/3$, the
number of edges between $\Gamma_{u,a}$ and $V_i$ is at least
$\left(\frac{\epsilon}{16} - \epsilon^3\right)|\Gamma_{u,a}| \cdot
|V_i|$. Before proceeding, we make the following useful
observation.

\begin{obs}\label{obs:bip}
Let $H=(A,B)$ be a bipartite graph with $\gamma |A||B|$ edges. Then
$H$ contains a matching $M$ of size $\gamma \frac{|A||B|}{|A|+|B|}$.
\end{obs}
\begin{proof}
Consider the following process that creates $M$. Initially $M_0 =
\emptyset$. Then in step $i$, we pick an arbitrary edge $\{a,b\}
\in E(H)$, set $M_{i+1} = M_i \cup \{a,b\}$ and remove from $E(H)$
all the edges incident with either $a$ or $b$. Clearly, in each
step the number of removed edges is bounded by $|A| + |B|$, so the
process continues for at least $\frac{E(H)}{|A|+|B|} = \gamma
\frac{ |A||B|}{|A|+|B|}$ steps. Hence $|M| = |\bigcup_i M_i| \geq
\gamma \frac{ |A||B|}{|A|+|B|}$.
\end{proof}

Returning to the proof of Lemma \ref{lem:many}, by Observation
\ref{obs:bip} the size of a maximal set $W_u$ as above is at least
$$\left(\frac{\epsilon}{16} -
\epsilon^3\right)\frac{|\Gamma_{u,a}| |V_i|}{|\Gamma_{u,a}| +
|V_i|} \geq \left(\frac{\epsilon}{16} -
\epsilon^3\right)\frac{\Big( \epsilon n /(3k)\Big) (n/k)}{\epsilon
n /(3k) + n/k} \geq \frac{\epsilon^2}{64 k} n.$$

To prove that $W = W_u \cap W_a$ is large, we similarly use the
regularity condition, but now on the pair $(\Gamma_{v,b},W_u)$. We
get,

$$|E(\Gamma_{v,b},W_u)| \geq \left(\frac{\epsilon}{16} -
\epsilon^3\right)|\Gamma_{v,b}| |W_u|.$$ Here too, by Observation
\ref{obs:bip} we can bound from below the size of a maximal
matching in the pair $(\Gamma_{v,b}, W_u)$ with
$$ \left(\frac{\epsilon}{16} -
\epsilon^3\right)\frac{|\Gamma_{v,b}| |W_u|}{|\Gamma_{v,b}| +
|W_u|} \geq \left(\frac{\epsilon}{16} -
\epsilon^3\right)\frac{\Big( \frac{\epsilon}{3k}n\Big)
\Big(\frac{\epsilon^2}{64 k} n\Big)}{\frac{\epsilon}{3k}n +
\frac{\epsilon^2}{64 k} n} \geq \epsilon^4\frac{n}{k} \geq
\epsilon^4\frac{n}{L} > 8^5 \log N,$$ where the last inequality
follows from our choice of $N$. Recall that the matching that we
found defines the desired set $W$, concluding the proof.
An algorithmic version of this lemma can be derived by simply using
an algorithmic version of Lemma \ref{lem:deg_new} in the selection
of $V_1,\ldots,V_k$ above. $\qed$

\subsection{Graphs with diameter $2$}
{\bf Proof of Theorem \ref{t-8logn}.}\, Consider a random
$3$-coloring of $E$, where every edge is colored with one of three
possible colors uniformly and independently at random. It is
enough to prove that for all pairs $u,v \in V$ the probability that
they are not rainbow connected is at most $1/n^2$. Then the proof
follows by the union bound (cf. \cite{AlSp}).

Let us fix a pair $u,v \in V$, and bound from above the probability
that this pair is not rainbow connected. We know that both $\Gamma(u)$
and $\Gamma(v)$ (the neighborhoods of $u$ and $v$) contain at least $8 \log n$ vertices.

\begin{enumerate}
\item
If $\{u,v\} \in E$ then we are done.
\item
If $|\Gamma(u) \cap \Gamma(v)| \geq 2 \log n$ then there are at
least $2 \log n$ edge-disjoint paths of length two from $u$ to $v$.
In this case, the probability that none of these paths is a rainbow
path is bounded by $(1/3)^{2 \log n} < 1/n^2$, and we are done.
\item
Otherwise, let $A = \Gamma(u) \setminus \Gamma(v)$ and $B = \Gamma(v) \setminus \Gamma(u)$.
We know that $|A|,|B| \geq 6 \log n$, and in addition, since the
first two cases do not hold and the diameter of $G$ is two, all the
(length two) shortest paths from $A$'s vertices to $v$ go through the
vertices in $B$. This implies that every vertex $x \in A$ has a
neighbor $b(x) \in B$ ($b(x)$ need not be a one-one function).
Let us consider the set of at least $6 \log
n$ edge-disjoint paths $P = \{u,x,b(x) : x \in A\}$. For each $x \in
A$, the probability that $u,x,b(x),v$ is a rainbow path (given the
color of the edge $(b(x),v)$) is $2/9$. Moreover, this event is
independent of the corresponding events for all other members of $A$,
because this probability does not change even with full knowledge
of the colors of all edges incident with $v$. Therefore, the probability that none of the
paths in $P$
extends to a rainbow path from $u$ to $v$ is at most $(7/9)^{6 \log
n} \leq 1/n^2$, as required.
\end{enumerate}

The above proof immediately implies a probabilistic polynomial expected time randomized algorithm
with zero error probability (since we can also efficiently check if the coloring indeed makes $G$
$3$-rainbow connected). The algorithm can be derandomized and converted to
a polynomial time probabilistic algorithm using the method of conditional expectations
(cf. \cite{AlSp}) similarly to the proof of Lemma \ref{lem:ab-rainbow}: For every partial coloring
of the edges we can efficiently bound the conditional probability that
a fixed pair $u,v$ is not rainbow-connected, using the relevant one
of the three cases concerning $u$ and $v$ that were analyzed above.
Now we can color the edges one by one, at each time taking care not to increase
the bound on the conditional expectation of unconnected pairs that results
from the above probability bound for every $u$ and $v$. Since the bound
on the expectation was smaller than $1$ before the beginning of the process,
in the end we would get a valid $3$-rainbow-coloring of $G$.
\qed

\section{Concluding remarks and open problems} \label{sec:concl}
\begin{itemize}
\item
Theorem \ref{thm:main} asserts that a connected graph with minimum
degree at least $\epsilon n$ has bounded rainbow connection.
However, the bound obtained is huge as it follows from the
Regularity Lemma. It would be interesting to find the ``correct''
bound. It is even possible that $rc(G) \le C/\epsilon$ for some
absolute constant $C$.
\item
The proof of Theorem \ref{t-rc2npc} shows that deciding whether $rc(G)=2$ is NP-Complete.
Although this suffices to deduce that computing $rc(G)$ is NP-Hard,
we still do not have a proof that deciding whether $rc(G) \le k$ is NP-Complete for every
fixed $k$. We can easily it for every {\em even} $k$ by the following reduction from the case
$k=2$. Given a graph $G$, subdivide every edge into $k/2$ edges.
Now, the new graph $G'$ has $rc(G')=k$ if and only if $rc(G)=2$.
Indeed, if $rc(G)=2$ then take a corresponding red-blue coloring of $G$ and color
$G'$ by coloring every subdivided red edge of $G$ with the colors $1,\ldots,k/2$ and every
subdivided blue edge with the colors $k/2+1,\ldots,k$. Conversely, if $G'$ has an edge
coloring making it rainbow connected using the colors $1,\ldots,k$, then color each edge
$e$ of $G$ as follows. If the subdivision of $e$ contains the color $1$, color $e$ red;
otherwise, color $e$ blue. This red-blue coloring of $G$ makes $G$ rainbow connected.

It is tempting to conjecture that for every $k$ it is NP-Hard even to
distinguish between $2$-rainbow-colorable graphs and graphs that are not
even $k$-rainbow-colorable.
\item
A parameter related to rainbow connection is the {\em rainbow
diameter}. In this case we ask for an edge coloring so that for any
two vertices, there is a rainbow {\em shortest} path connecting
them. The rainbow diameter, denoted $rd(G)$ is the smallest
number of colors used in such a coloring. Clearly, $rd(G) \ge rc(G)$
and obviously every connected graph with $n$ vertices has $rd(G) <
\binom{n}{2}$. Unlike rainbow connection, which is a monotone
graph property (adding edges never increases the rainbow connection)
this is not the case for the rainbow diameter (although we note that
constructing an example that proves non-monotonicity is not
straightforward). Clearly, computing $rd(G)$ is NP-Hard since
$rc(G)=2$ if and only if $rd(G)=2$. It would be interesting to prove
a version of Theorem \ref{thm:main} for rainbow diameter. We
conjecture that, indeed, if $G$ is a connected graph with minimum
degree at least $\epsilon n$ then it has a bounded rainbow diameter.
\item
Suppose that we are given a graph $G$ for which we are {\em told} that
$rc(G)=2$. Can we rainbow-color it in polynomial time with $o(n)$
colors? For the usual coloring problem, this version has been well
studied. It is known that if a graph is $3$-colorable (in the usual
sense), then there is a polynomial time algorithm that colors it with
$\tilde{O}(n^{3/14})$ colors \cite{BlKa}.
\end{itemize}

\end{document}